\documentclass[a4paper,11pt]{article}
\usepackage[multidot]{grffile}
\usepackage{latexsym,amssymb,enumerate,amsmath,epsfig,amsthm}
\usepackage[margin=1in]{geometry}
\usepackage{setspace,color}
\usepackage{graphics}
\usepackage{graphicx}
\usepackage[ruled]{algorithm2e}
\usepackage{empheq}
\usepackage{bm,amsmath}
\usepackage{subcaption}
\usepackage{hyperref} 
\hypersetup{
    colorlinks=true,
    linkcolor=blue,
    filecolor=magenta,      
    urlcolor=cyan,
    citecolor=red,
    pdfpagemode=FullScreen,
    }

\usepackage{rotating,booktabs,array,amsmath,amssymb}
\newcolumntype{C}{>{$\displaystyle}c<{$}} 

\usepackage{multirow}
\newcommand\fixup{\kern-\fontcharic\scriptfont2`\"}

\usepackage{lineno}

\newcommand{\x}{\mathbf{x}}
\newcommand{\s}{\mathbf{s}}
\newcommand{\p}{\mathbf{p}}
\newcommand{\q}{\mathbf{q}}

\title{A Spherical Crank-Nicolson Integrator Based on the Exponential Map and the Spherical Linear Interpolation}

\author{
Shingyu Leung\thanks{Department of Mathematics, the Hong Kong University of Science and Technology, Clear Water Bay, Hong Kong. Email: {\bf masyleung@ust.hk}}
}

\date{}

\begin{document}
\thispagestyle{plain}
\maketitle

\begin{abstract}
We propose implicit integrators for solving stiff differential equations on unit spheres. Our approach extends the standard backward Euler and Crank-Nicolson methods in Cartesian space by incorporating the geometric constraint inherent to the unit sphere without additional projection steps to enforce the unit length constraint on the solution. We construct these algorithms using the exponential map and spherical linear interpolation (SLERP) formula on the unit sphere. Specifically, we introduce a spherical backward Euler method, a projected backward Euler method, and a second-order symplectic spherical Crank-Nicolson method. While all methods require solving a system of nonlinear equations to advance the solution to the next time step, these nonlinear systems can be efficiently solved using Newton's iterations. We will present several numerical examples to demonstrate the effectiveness and convergence of these numerical schemes. These examples will illustrate the advantages of our proposed methods in accurately capturing the dynamics of stiff systems on unit spheres.
\end{abstract}

\section{Introduction}

We consider solving the ordinary differential equation (ODE) on $\mathbb{S}^2$ given by $\p'(t)=f(\p(t),t)$ with an initial condition $\p(0) = \p_0 \in \mathbb{S}^2$. Here, $f:\mathbb{S}^2\times [0,\infty) \rightarrow \mathcal{T}(\mathbb{S}^2)$ satisfies a Lipschitz function ensuring that the solution to the ODE remains confined to $\mathbb{S}^2$ for all time. This formulation finds natural application in scenarios such as path planning for rigid bodies, where the ODE of the special orthogonal group $SO(3)$ needs to be solved \cite{shi09}. Other applications include quantum field theory within quantum mechanics \cite{adl86}, protein structure modeling \cite{pro14}, molecular dynamics simulation \cite{rap85}, fluid mechanics theory \cite{ghkr06}, fluid flow visualization \cite{hanma95}, computations involving flexible filaments and fibers in complex fluids \cite{tscfro20,stwk21}, differential equations \cite{kouxia18}, and dynamics of rigid bodies \cite{weiterfed06,wil09,udwsch10}. Solutions to this problem will also be useful in various applications, particularly high-frequency wave propagation on a unit sphere utilizing geometrical optics \cite{gri68}, spin dynamics based on the Landau-Lifshitz equation \cite{frahualei97,krebunlan98,ewan00,ban05,mclone06,cim07,jeokim10,ser14,ccwx22}, and $p$-harmonic flows for signal denoising \cite{tansapcas00,tansapcas01,vesosh02,lysoshtai04,golwenyin09}. 

Equations in this form have been of great importance, leading to the proposal of various numerical approaches. When solving the Landau-Lifshitz equation directly in $\mathbb{R}^3$, \cite{frahualei97} introduced multiple low-order geometric integrators by leveraging the Hamiltonian structure of the system. Geometric integrators based on Lie group methods have been developed in \cite{lewnig03}. Another approach, presented in \cite{ewan00}, involved the development of a simple projection method. The work by \cite{ban05} introduced a projected backward Euler method. Furthermore, \cite{krebunlan98} suggested a higher-order predictor-corrector approach. A recent approach based on operator splitting with projection has been proposed in \cite{ccwx22}. Recognizing the importance of preserving the total energy in Hamiltonian systems, efforts have also been made to construct second-order midpoint-type implicit schemes for solving equations on the unit sphere. In this regard, \cite{jeokim10} developed a direct extension of the Crank-Nicolson method from $\mathbb{R}^3$. This method was demonstrated to satisfy the geometric constraint for specific velocities.

Several more recent approaches have been proposed in the literature. One such approach is a symmetric integrator developed in \cite{mclmodver14,mclmodver17,dzmbbb22}, given by the equation
$$
\p^{n+1} = \p^{n} + h f\left(\frac{\p^{n+1}+\p^{n}}{\|\p^{n+1}+\p^{n}\|}, t^{n+\frac{1}{2}} \right) \, .
$$
To evaluate the velocity, this method first takes the average of the two intermediate solutions. However, since the \textit{midpoint} of these two points is generally not on the unit sphere, the method incorporates an extra projection step to evaluate the associated velocity. The main issue with this approach is that the update formula is essentially a simple mimic of the Euler method in $\mathbb{R}^3$. In general, the velocity at the midpoint is not on the tangent plane at either $\mathbf{p}^n$ or $\mathbf{p}^{n+1}$. As a result, the computed $\mathbf{p}^{n}$ using this formula may not lie on the unit sphere. This is why a further projection step is necessary for the computed $\mathbf{p}^{n+1}$.

Another suggested approach is based on the backward differentiation formula (BDF) \cite{xgwzc20}. Using the second-order BDF formula, we have
$$
\frac{3}{2}\p^{n+2} -2\p^{n+1}+\frac{1}{2} \p^n = h f\left(\p^{n+2}, t^{n+2} \right) \, .
$$
It is important to note that the BDF expression on the left-hand side of the scheme is derived for the space $\mathbb{R}^3$. Therefore, there is no guarantee that the linear combination of these approximated points will give a vector lying on the tangent plane at $\mathbf{p}^{n+2}$. The resulting vector may not be properly aligned with the geometry of the unit sphere. An additional projection step is, therefore, necessary to ensure that the computed vector lies on the tangent plane and satisfies the unit sphere constraint.

Designed for general manifolds, \cite{hai00} has developed an interesting symmetric projection approach that can retain the time-reversibility property. The idea is to allow intermediate solutions to leave the manifold while adjusting the perturbation simultaneously to make the overall procedure symmetric. The work by \cite{ceor20} has developed high-order energy-preserving collocation-like methods on Riemannian manifolds.

In a recent paper \cite{leuchalee24}, we introduce a class of explicit Runge-Kutta (RK) integrators designed to solve ODEs on the unit sphere up to third-order accuracy. Unlike the conventional projected RK methods that require an additional projection step, our approach leverages the explicit formula of the exponential map on the unit sphere. This approach yields a numerical solution that automatically satisfies the $\mathbb{S}^2$ constraint. We adopt the TVDRK methodology \cite{shu88,shuosh88,gotshu98,gotshutad00}, which constructs higher-order numerical solutions using convex combinations of elementary forward Euler-type building blocks. This construction gives rise to a class of straightforward, high-order, and efficient explicit numerical schemes for solving ODEs on $\mathbb{S}^2$. 

In certain applications, however, implicit integrators are preferred due to their superior stability properties. Consequently, in this work, we propose implicit numerical integrators that rely on the exponential map and the spherical linear interpolation (SLERP) formula of unit spheres. These methods ensure solutions automatically satisfy the geometrical constraint without any explicit projection step. In particular, we introduce two first-order methods inspired by the backward Euler method, as well as a second-order Crank-Nicolson method based on the midpoint approach. Similar to conventional Cartesian implicit methods, our proposed approaches involve solving a system of nonlinear equations. To address this, we develop iterative methods utilizing Newton's iteration, which offers a numerically efficient approach. 

The rest of the paper is structured as follows. In Section \ref{Sec:Background}, we will begin by summarizing the SLERP formula, which will help express our proposed spherical Crank-Nicolson scheme. Additionally, we will present the spherical forward Euler method developed in \cite{leuchalee24}, which will serve as a basis for comparing the performance of the numerical integrators. In Section \ref{Sec:Proposed}, we will introduce three proposed implicit integrators on spheres. Finally, in Section \ref{Sec:Examples}, we will present several numerical examples to demonstrate the effectiveness and convergence of our numerical algorithms.


\section{Background}
\label{Sec:Background}

\subsection{The Spherical Linear Interpolation (SLERP) Formula}
\label{SubSec:SLERP}

This section provides the background on the interpolation of spherical data and introduces the spherical linear interpolation (SLERP) formula \cite{shoemake_85,sola} using a quaternion representation \cite{ham63}. Quaternions are numbers consisting of four dimensions, one real part, and a three-dimensional analogy to {the imaginary part of} complex numbers. A quaternion can be written in many forms: 
$
\underset{\text{real}}{\boxed{a}} + \underset{\text{imaginary}}{\boxed{b \mathbf{i} + c\mathbf{j} + d\mathbf{k}}} = (a, b, c, d) = (\underset{\text{scalar}}{\boxed{a}}, \underset{\text{vector}}{\boxed{\mathbf{u}}}) \, , 
$ where \(a, b, c, d \in \mathbb{R}\), $\mathbf{u} = (b, c, d) \in \mathbb{R}^3$. The notations $\mathbf{i}$, $\mathbf{j}$, and $\mathbf{k}$ are extensions of {the imaginary part of} complex numbers. Based on the quaternion representation, the SLERP (\textit{Spherical Linear intERPolation}) formula can be expressed by
$
\mbox{SLERP}(\mathbf{\mathbf{q_a}}, \mathbf{q_b}, t) = (\mathbf{q_a}) ((\mathbf{q_a})^{-1} \mathbf{q_b})^t 
$
where we have applied the following quaternion properties
\begin{itemize}
\item Hamilton product: $(a_1, \mathbf{u_1})(a_2, \mathbf{u_2}) = (a_1a_2 - \mathbf{u_1} \cdot \mathbf{u_2}, a_1\mathbf{u_2} + a_2\mathbf{u_1} + \mathbf{u_1} \times \mathbf{u_2})$ where the notation {$\cdot$ and $\times$ denotes the typical dot and cross product.}
\item Inverse map: $\mathbf{q}^{-1} = (a, -\mathbf{u})/(a^2 + b^2 + c^2 + d^2)$.
\item Exponential map: $\exp(a, \mathbf{u}) = \exp(a)(\cos \lVert\mathbf{u}\rVert, ((\sin \lVert\mathbf{u}\rVert)/\lVert\mathbf{u}\rVert) \mathbf{u})$ where $\lVert \cdot \rVert$ denotes the 2-norm.
\item Logarithm map: $ \ln(a, \mathbf{u}) = \left(\ln \sqrt{a^2 + \lVert \mathbf{u} \rVert^2}, \frac{1}{\lVert \mathbf{u} \rVert} 
\arccos \left( \frac{a}{\sqrt{a^2 + \lVert \mathbf{u} \rVert^2}} \right) \mathbf{u} \right)$. 
\item Power map: $(a, \mathbf{u})^{f(t)} = \exp(f(t)\ln(a, \mathbf{u}))$.
\end{itemize}

\subsection{Spherical Forward Euler (SFE)}

Consider moving the data point $\p^n$ with a nonzero constant velocity $\s = f(\p^n,t^n)$ for the period $h$. The arrival location on the unit sphere has an explicit formula given by $\p^{n+1} = \mbox{exp}_{\p^{n}}(h\s)$ where $\mbox{exp}_{\p}: T_{\p}\mathbb{S}^2 \rightarrow \mathbb{S}^2$ is the exponential map with $\mbox{exp}_{\p}(\s)=\gamma(1)$ where $\gamma$ being the unique geodesic satisfying $\gamma(0)=\p$ and $\gamma'(0)=\s$. Mathematically, we have the following expression explicitly for the unit sphere $\mbox{exp}_{\p}(\s) = \cos(\|\s\|) \p + \sin(\|\s\|) \frac{\s}{\|\s\|}$, and therefore, we arrive the spherical forward Euler (SFE) method
$$
\boxed{
\mbox{Spherical FE (SFE): } \, \left\{
\begin{array}{l}
\s_1 = f(\p^n,t^n) \\
\p^{n+1} = \cos(h \|\s_1\|) \p^n + \sin(h\|\s_1\|) \frac{\s_1}{\|\s_1\|} \, .
\end{array}
\right.
}
$$
The SFE scheme emulates the conventional forward Euler method utilized for solving ODEs in Cartesian space. However, in contrast to permitting the solution to reach any arbitrary point in the entire space after a single timestep, the SFE method respects the spherical geometry by ensuring that the solution remains confined to $\mathbb{S}^2$. We can easily prove using Taylor series expansion that this scheme possesses local second-order accuracy, resulting in a globally first-order accurate solution. A more detailed proof can be found in \cite{leuchalee24}.

\section{Our Proposed Implicit Integrators}
\label{Sec:Proposed}

\subsection{A Spherical Backward Euler Method}

In this section, we will construct a spherical backward Euler method. The idea is to determine $\p^{n+1}$ implicitly by solving the following system of equations:
\begin{equation}
 \boxed{
\mbox{Spherical BE (SBE): } \, 
 \left\{
\begin{array}{l}
\s = f(\p^{n+1},t^{n+1}) \\
\p^n = \mbox{exp}_{\p^{n+1}}(-h \s) = \cos(h \|\s\|) \p^{n+1} - \sin(h\|\s\|) \frac{\s}{\|\s\|} \, .
\end{array}
\right.
}
\label{Eqn:SBE}
\end{equation}
The first equation evaluates the velocity defined on the tangent plane at the unknown location $\p^{n+1}$, while the second equation indicates that when tracing the exponential map \textit{backward} in time from $\p^{n+1}$ with the specific terminal velocity $\s$, one will reach the takeoff location $\p^n$ after a time period of $h$.

A simple approach is to solve this system of nonlinear equations using an alternating iterative approach. We note that the system (\ref{Eqn:SBE}) is already in a form that is best suited for a fixed-point iterative approach, given by
$$
 \left\{
\begin{array}{l}
\s_{k+1} = f(\q_{k},t^{n+1}) \\
\q_* = \mbox{exp}_{\q_{k+1}}(-h \s_{k+1}) 
\end{array}
\right.
$$
with the index $k \ge 0$ and the point $\q_* = \p^n$. The second stage of this iteration is, in fact, linear in $\q_{k+1}$. We can actually solve the exponential map and obtain
$$
\q_{k+1} = \frac{1}{\cos(h \|\s_{k+1}\|)} \left[ \q_*+ \sin (h \|\s_{k+1}\|) \frac{\s_{k+1}}{\|\s_{k+1}\|} \right] \, .
$$
Once this iteration converges, such that, for example, $|\q_{k+1}-\q_k|<\epsilon$, we assign the approximation to the differential equation at $t=t^{n+1}$ as $\p^{n+1}=\q_k$. There are two main issues with this approach. The first issue concerns the convergence of the iterative approach. Since the iterative function is not straightforward, there is no guarantee of convergence in general. The second issue concerns the constraint that the numerical solution $\p^{n+1}$ should remain on the unit sphere. Specifically, there is no guarantee that any intermediate step will yield $\q_{k+1}$ that stays on the unit sphere, even if $\q_k$ does. Examining the update formula for $\q_{k+1}$, we observe that the velocity vector $\s_{k+1}$ lies on the tangent plane at $\q_k$, rather than $\q_*$. Unless the point $\q_k$ solves the system (\ref{Eqn:SBE}), such a velocity will generally pull $\q_*$ away from the unit sphere. When $\q_{k+1}$ does not have unit length, it may be challenging to evaluate $f(\q_{k+1},t^{n+1})$ in the subsequent fixed-point iteration. This can be easily corrected by incorporating an additional projection step, resulting in
$$
 \left\{
\begin{array}{l}
\s_{k+1} = f(\q_{k},t^{n+1}) \\
\q_{k+1/2} = \frac{1}{\cos(h \|\s_{k+1}\|)} \left[ \q_*+ \sin (h \|\s_{k+1}\|) \frac{\s_{k+1}}{\|\s_{k+1}\|} \right] \, , \, \q_{k+1}= \frac{\q_{k+1/2}}{\|\q_{k+1/2}\|} \, .
\end{array}
\right.
$$
However, it is important to note that incorporating this extra projection step may introduce additional complications regarding the convergence of the entire fixed-point iteration.

Instead, we propose to extend the nonlinear system to $\mathbb{R}^6$ directly and solve the equations using Newton's iterations at once. The implementation of this approach is straightforward, and the convergence of the iteration is extremely fast. We rewrite equation (\ref{Eqn:SBE}) as
\begin{equation}
 \left\{
\begin{array}{l}
\s - g(\q) =0 \\
-\q_* + \cos(h \|\s\|) \q - \sin(h\|\s\|) \frac{\s}{\|\s\|} =0 
\end{array}
\right.
\label{Eqn:SBE2}
\end{equation}
where we define $g(\q)=f(\q,t^{n+1})$ for simplicity. The Jacobian $D(\s,\q)\in\mathbb{R}^{6\times 6}$ of this nonlinear system is given by
$$
D(\s,\q) = 
\left(
\begin{array}{c|c}
I_3 & -G(\q) \\ \hline
J(\s,\q) & \cos(h\|\s\|) I_3
\end{array}
\right) \, .
$$
Here, $I_3$ is the $3 \times 3$ identity matrix, the individual elements of the matrix $G$ are given by $G_{i,j} = \partial g_i / \partial q_j$ with $g(\q) = (g_1(\q), g_2(\q), g_3(\q))$ and $\q = (q_1, q_2, q_3)$ for $i,j = 1,2,3$. The elements of $J(\s,\q)$ are given by
$$
J_{i,j}=-h \sin(h\|\s\|) \frac{q_i s_j}{\|\s\|} - h \cos(h \|\s\|) \frac{s_is_j}{\|\s\|^2} - \sin (h \|\s\|) \left( \frac{\delta_{i,j}}{\|\s\|} - \frac{s_i s_j}{\|\s\|^3} \right) 
$$
where $\delta_{i,j} = 1$ when $i$ equals $j$ and is zero otherwise. Instead of defining an extension of $g(\q)$ for general $\q\in\mathbb{R}^3$, we suggest incorporating a projection step to ensure that the intermediate solution $\q_k$ remains on the unit sphere. Therefore, the projected Newton's iteration is given by
$$
 \left\{
\begin{array}{l} 
\left(
\begin{array}{c}
\s_{k+1} \\
\q_{k+1/2}
\end{array}
\right)
 = 
 \left(
\begin{array}{c}
\s_{k} \\
\q_{k}
\end{array}
\right)
- D(\s_k,\q_k)^{-1}
 \left[
\begin{array}{c}
\s_{k} - g(\q_{k}) \\
-\q_* + \cos(h \|\s_k\|) \q_k - \sin(h\|\s_k\|) \frac{\s_k}{\|\s_k\|}
\end{array}
\right] \\
\q_{k+1}=\frac{\q_{k+1/2}}{\|\q_{k+1/2}\|} \, .
\end{array}
\right.
$$
We assign the spherical forward Euler solution and the velocity at the corresponding location as the initial guess of the iteration. Specifically, we set $\q_0=\exp_{\p^n}(h\|f(\q_*,t^n)\|)$ and $\s_0=g(\q_0)$. To avoid issues of division by zero, we can modify the computation of the norm of the velocity $\s_k$ by replacing it with $\max(\|\s_k\|, \epsilon_{\mbox{\tiny machine}})$, where $\epsilon_{\mbox{\tiny machine}}$ is the machine epsilon. This modification ensures that the computation remains stable even at stationary points on the sphere.

\begin{figure}
\centering
\includegraphics[width=0.99\textwidth]{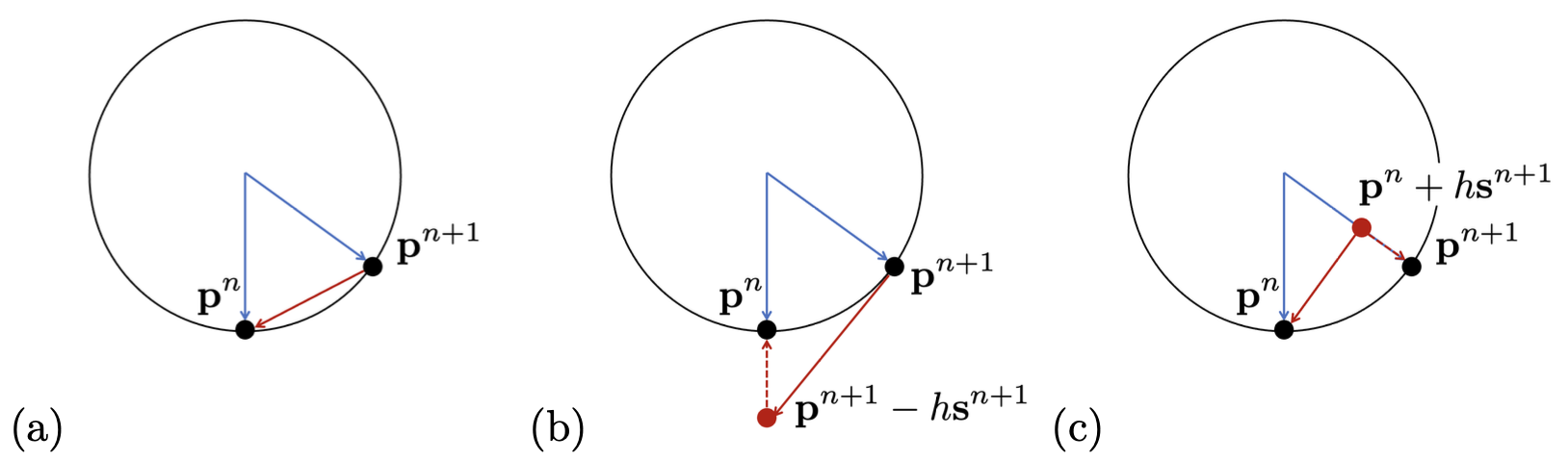}
\caption{(Section \ref{SubSec:BackwardEuler})(a) A simple backward Euler approach (\ref{Eqn:BackwardEuler1}) does not work. (b) A slightly more involved backward Euler implementation (\ref{Eqn:BackwardEuler2}). (c) Our proposed projected backward Euler method (\ref{Eqn:BackwardEuler3}).} 
\label{Fig:BackwardEuler}
\end{figure}

\subsection{Some Projected Backward Euler Methods}
\label{SubSec:BackwardEuler}

It is also possible to design backward Euler methods without using the exponential map. Although they are not the main focus of this paper, the idea might be helpful when dealing with other manifolds.

The first idea is to determine the location $\p^{n+1} \in \mathbb{S}^2$ such that a backward Euler step reaches a point that collides with $\p^n$ \textit{without} using any projection, as shown in Figure \ref{Fig:BackwardEuler}(a). Mathematically, this implies finding a location $\q \in \mathbb{S}^2$ that satisfies the following system of equations:
\begin{equation}
 \left\{
\begin{array}{l}
\s = g(\q) \, , \, \p^n = \q-h \s \, .
\end{array}
\right.
\label{Eqn:BackwardEuler1}
\end{equation}
However, the system has no solution on the unit sphere for $h > 0$ and a general velocity field $\|\s\|>0$. This is because $\q \cdot \s=0$ and $\|\p^n\|^2 = \|\q - h\s\|^2 = \|\q\|^2 + h^2\|\s\|^2 = 1 + h^2\|\s\|^2 > 1$, which implies that $\p^n$ cannot have a unit length.

Another possible scheme is to determine a $\p^{n+1}\in\mathbb{S}^2$ such that a backward Euler step reaches a point that collides with $\p^n$ \textit{after} a projection step, as shown in Figure \ref{Fig:BackwardEuler}(b). Mathematically, we solve the following nonlinear system
\begin{equation}
 \left\{
\begin{array}{l}
\s = g(\q) \, , \, \p^n = \left(\q-h \s\right)/\|\q-h \s\| 
\end{array}
\right.
\label{Eqn:BackwardEuler2}
\end{equation}
for $(\s,\q)$. Once a solution $(\s^*,\q^*)$ is determined, we assign $\p^{n+1}=\q^*$. This approach, however, leads to a complicated nonlinear equation, especially when solving by Newton's iterations. The computation of the derivative of the right-hand side of the second equation with respect to $\q$ (and also with respect to $\s$ if one embeds the approach to $\mathbb{R}^6$) can make the iterations more involved.

Instead, we consider the third interpretation that we first extend the velocity into the whole space and look for a point in $\mathbb{R}^3$ so that a backward Euler step will reach $\p^n$. Then, we project this point back to the unit sphere and assign the projected point as $\p^{n+1}$, as shown in Figure \ref{Fig:BackwardEuler}(c). Mathematically, this implies solving
\begin{equation}
\boxed{
\mbox{Projected BE (PBE): } \, 
 \left\{
\begin{array}{l}
\s = g(\q/\|\q\|) \\
\p^n = \q-h \s
\end{array}
\right.
}
\label{Eqn:BackwardEuler3}
\end{equation}
Once a solution $(\s^*,\q^*)$ is determined, we assign $\p^{n+1}=\q^*/\|\q^*\|$. Since the expression is rather simple, one might implement a Newton's iteration on $\q$ instead of $(\s,\q)$. In particular, we obtain the following iterative form
$$
\q_{k+1}=\q_k-D(\q_k)^{-1} \left[-\q_*+\q_k-hg(\q^k)\right]
$$
where $D(\q)=I_3-hG(\q)$. Investigating this iterative formula, we obvious that when $h$ is large, the matrix $D$ might be closed to singular which might be undesired in some applications. Therefore, it also becomes clear the advantages of embedding the nonlinear system in $\mathbb{R}^6$ instead. This gives the following iterative form,
$$
 \left\{
\begin{array}{l} 
\left(
\begin{array}{c}
\s_{k+1} \\
\q_{k+1/2}
\end{array}
\right)
 = 
 \left(
\begin{array}{c}
\s_{k} \\
\q_{k}
\end{array}
\right)
- D(\s_k,\q_k)^{-1}
 \left[
\begin{array}{c}
\s_{k} - g(\q_{k}) \\
-\q_*+\q_k-hg(\q^k)
\end{array}
\right] \\
\q_{k+1}=P(\q_{k+1/2})
\end{array}
\right.
$$
where $P(\q)$ is the projection of $\q$ onto the unit sphere and the Jacobian is given by
$$
D(\s,\q) = 
\left(
\begin{array}{c|c}
I_3 & -G(\q) \\ \hline
-hI_3 & I_3
\end{array}
\right) \, .
$$

As a final remark on the typical backward Euler method, it should be noted that this approach can also be applied to differential equations on general manifolds. The only modification required is to redefine the projection operation $P(\q)$. Instead of using $P(\q) = \q/\|\q\|$ as in the case of the unit sphere, we can determine the closest point projection onto the manifold $\Sigma$ by assigning $P(\q) = \text{argmin}_{\p \in \Sigma} \|\q - \p\|$ or $P(\q) =  \text{argmin}_{\p \in \Sigma} \frac{1}{2}\|\q - \p\|^2$.

\begin{figure}
\centering
\includegraphics[width=0.99\textwidth]{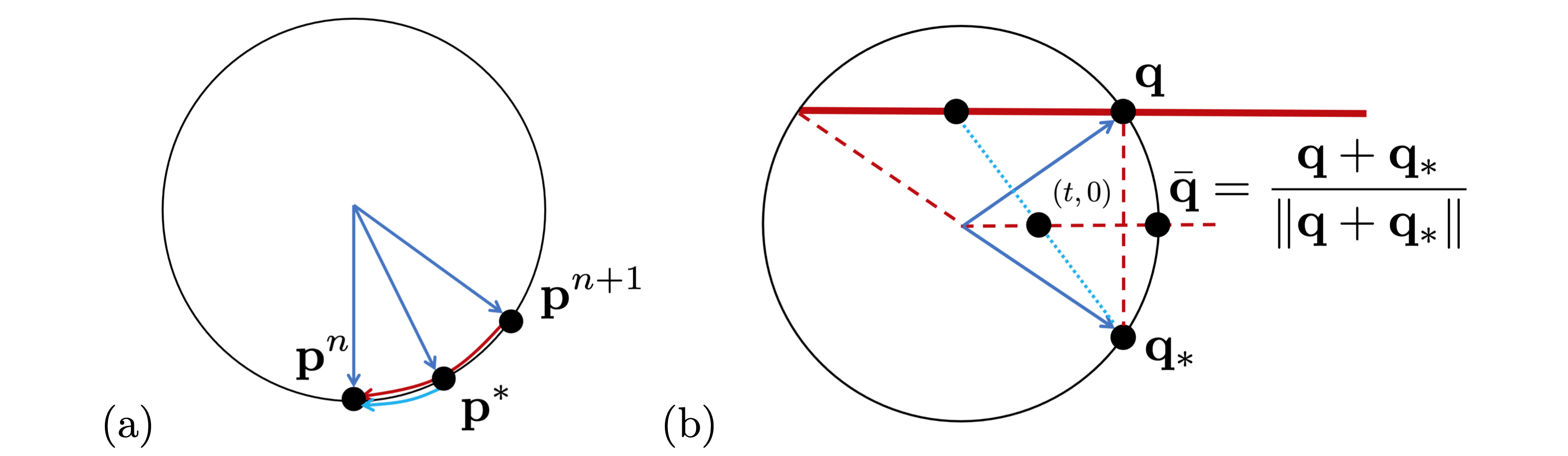}
\caption{(Section \ref{SubSec:SphericalCN}) (a) The spherical Crank-Nicolson scheme. (b) The solution to (\ref{Eqn:CNNotWork}) is not unique.} 
\label{Fig:SphericalCN}
\end{figure}

\subsection{A Spherical Crank-Nicolson Method}
\label{SubSec:SphericalCN}

This section develops a new second-order fully implicit integrator on spheres based on the Crank-Nicolson or the midpoint-type idea. We propose the following scheme by solving the following system of nonlinear equations for both $(\s,\q)$
\begin{equation}
\boxed{
\mbox{Spherical CN (SCN): } \, 
 \left\{
\begin{array}{l}
\s = f \left(\p^{*},t^{n+\frac{1}{2}} \right) \, , \, \p^{*} = \mbox{SLERP}\left(\p^n,\q,\frac{1}{2}\right) \, , \\
\p^n = \mbox{exp}_{\p^{*}} \left( -\frac{h}{2} \s \right)  \, .
\end{array}
\right.
}
\label{Eqn:SCN}
\end{equation}
Once we have the solution $(\s^*,\q^*)$, we assign $\p^{n+1}=\q^*=(q_{*1},q_{*2},q_{*3})$. In the scheme, $\p^*$ represents the midpoint along the geodesic between $\p^n$ and the unknown location $\q$ (which corresponds to the numerical approximation of the evolution at $t=t^{n+1}$). The vector $\s$ denotes the velocity at this midpoint $\p^*$, which also lies on the tangent plane at the same point. The second equation in the system indicates that when tracing \textit{backward} in time for half a step size from this midpoint $\p^*$ with the corresponding velocity using the SBE method, we should reach the starting point at $\p^n$.

Now, overloading the notation of $g$, we introduce $g(\q)=f \left(\q,t^{n+\frac{1}{2}} \right)$ and obtain the following nonlinear system,
\begin{equation}
 \left\{
\begin{array}{l}
\s - g\left[ \mbox{SLERP}\left(\p^n,\q,\frac{1}{2}\right) \right] =0 \, , \\
- \q_* + \cos\left( \frac{h}{2} \|\s\|\right) \mbox{SLERP}\left(\p^n,\q,\frac{1}{2}\right) - \sin\left(\frac{h}{2} \|\s\| \right) \frac{\s}{\|\s\|} =0  \, .
\end{array}
\right.
\label{Eqn:SCN2}
\end{equation}
Compared to the backward Euler method (\ref{Eqn:BackwardEuler3}) and its corresponding spherical version (\ref{Eqn:SBE2}), these equations are similar, and we only need to make slight modifications to the previous implementation. In particular, we have the following update formula:
$$
\footnotesize	
 \left\{
\begin{array}{l} 
\left(
\begin{array}{c}
\s_{k+1} \\
\q_{k+1/2}
\end{array}
\right)
 = 
 \left(
\begin{array}{c}
\s_{k} \\
\q_{k}
\end{array}
\right)
- D(\s_k,\q_k;\q_*)^{-1}
 \left[
\begin{array}{c}
\s_{k} - g \left[ \mbox{SLERP}\left(\q_*,\q_k,\frac{1}{2}\right) \right] \\
-\q_* + \cos \left( \frac{h}{2} \|\s_k\|\right) \mbox{SLERP}\left(\q_*,\q_k,\frac{1}{2}\right) - \sin \left(\frac{h}{2}\|\s_k\|\right) \frac{\s_k}{\|\s_k\|}
\end{array}
\right] \\
\q_{k+1}=\frac{\q_{k+1/2}}{\|\q_{k+1/2}\|} 
\end{array}
\right. \, .
$$
The Jacobian $D(\s,\q)$ given by
$$
D(\s,\q;\q_*) = 
\left(
\begin{array}{c|c}
I_3 & -G(\q;\q_*) \\ \hline
H(\s,\q;\q_*) & \cos \left( \frac{h}{2} \|\s\| \right) K(\q;\q_*)
\end{array}
\right) \, .
$$
with
\begin{eqnarray*}
G_{i,j}&=&\frac{\partial}{\partial q_j}  g_i \left[ \mbox{SLERP}\left(\q_*,\q,\frac{1}{2}\right)\right] \, , \\
H_{i,j}&=& - \frac{h}{2} \sin \left( \frac{h}{2}\|\s\| \right) \frac{s_j}{\|\s\|} \, \mbox{SLERP}_i\left(\q_*,\q,\frac{1}{2}\right)\\
& & - \frac{h}{2} \cos\left(\frac{h}{2} \|\s\|\right) \frac{s_is_j}{\|\s\|^2} - \sin \left(\frac{h}{2} \|\s\|\right) \left( \frac{\delta_{i,j}}{\|\s\|} - \frac{s_i s_j}{\|\s\|^3} \right) \, , \\
K_{i,j}&=&\frac{\partial}{\partial q_j} \left[ \mbox{SLERP}_i \left(\q_*,\q,\frac{1}{2}\right) \right] = \frac{1}{2} \sec \left( \frac{\theta}{2} \right) \left[ 1- \frac{1}{2}q_{*j}(q_{*i}+q_i) \sec^2 \left( \frac{\theta}{2} \right) \right]
\, ,
\end{eqnarray*}
where $\theta=\cos^{-1}(\q\cdot\q_*)$ is the angle between the two points on the unit sphere. To initiate the iteration, we use the SBE solution as the initial condition for the SCN iteration. We observe that the convergence of this Newton's iteration is extremely fast. In all of the numerical examples below, we find that it is sufficient to perform two to three iterations to achieve the necessary accuracy.

We observe that the SCN integrator is time-reversible. Specifically, the numerical solution $(\s^*,\p^n,\p^{n+1})$ satisfies both systems
$$
 \left\{
\begin{array}{l}
\s^* = f \left(\p^{*},t^{n+\frac{1}{2}} \right) \, , \, \p^{*} = \mbox{SLERP}\left(\p^n,\p^{n+1},\frac{1}{2}\right) \, , \, \p^n = \mbox{exp}_{\p^{*}} \left( -\frac{h}{2} \s \right)  
\end{array}
\right.
$$
and
$$
 \left\{
\begin{array}{l}
\s^* = f \left(\p^{*},t^{n+\frac{1}{2}} \right) \, , \, \p^{*} = \mbox{SLERP}\left(\p^{n+1},\p^{n},\frac{1}{2}\right) \, , \, \p^{n+1} = \mbox{exp}_{\p^{*}} \left( \frac{h}{2} \s \right)  \, ,
\end{array}
\right.
$$
since $\mbox{SLERP}\left(\mathbf{a},\mathbf{b},\alpha \right) = \mbox{SLERP}\left(\mathbf{b},\mathbf{a},1-\alpha \right)$ and the point $\p^*$ is the midpoint along the geodesic joining $\p^n$ and $\p^{n+1}$ with the linearized velocity $\s$. This property is important as it ensures the long-term stability and accuracy of the numerical solution, making it a reliable method for simulating systems with conserved quantities such as Hamiltonian systems. As we will demonstrate in Section \ref{SubSec:Ham}, this numerical approach not only gives a second-order accurate solution but also preserves the Hamiltonian in certain systems well. 

The SCN scheme introduces an additional stability constraint compared to the SBE method due to the SLERP interpolation. Specifically, when the points $\mathbf{a}$ and $\mathbf{b}$ are antipodes on the unit sphere, the interpolation $\mbox{SLERP}\left(\mathbf{a},\mathbf{b},\frac{1}{2} \right)$ is not uniquely determined since it is not possible to find a unique geodesic on the unit sphere that connects these two antipodal points. In the current application, the midpoint given by
$$
\mbox{SLERP} \left( \mathbf{p}^n,\exp_{\mathbf{p}^n}(h\|\mathbf{s}\|),\frac{1}{2} \right) 
$$
does not coincide with $\exp_{\mathbf{p}^n}\left(\frac{h}{2}\|\mathbf{s}\| \right)$ when $h\|\mathbf{s}\|>\pi$. This is because the midpoint does not lie on the trajectory $\gamma(t)$ for $t\in[0,1]$ satisfying $\gamma(0)=\mathbf{p}^n$ and $\gamma'(0)=\mathbf{s}$ when $\gamma(t)$ passes through the antipole for some $t\in[0,1]$. Therefore, although the conventional CN method is unconditionally stable for an arbitrary step size, we do not recommend choosing an extremely large step size for the SCN scheme.

To conclude this section, we mention that there could be various possible variations of the SCN scheme. For example, instead of connecting the starting point $\p^n$ with the midpoint along the geodesic $\p^*$ through the second equation in the nonlinear system (\ref{Eqn:SCN}), one could also replace it with the second half of the relation, which states that the trajectory will reach the unknown location $\q$ in half of a time step from the midpoint location $\p^*$ using the corresponding velocity $\s$. Mathematically, this can be expressed as $\q = \mbox{exp}_{\p^{*}} \left( \frac{h}{2} \s \right)$ which implies the following nonlinear system,
$$
 \left\{
\begin{array}{l}
\s - g\left[ \mbox{SLERP}\left(\p^n,\q,\frac{1}{2}\right) \right] =0 \, , \, \q - \cos\left( \frac{h}{2} \|\s\|\right) \mbox{SLERP}\left(\p^n,\q,\frac{1}{2}\right) - \sin\left(\frac{h}{2} \|\s\| \right) \frac{\s}{\|\s\|} =0  \, .
\end{array}
\right.
$$
This nonlinear system is similar to the one we developed in (\ref{Eqn:SCN2}) except that we replace the \textit{constant} term $-\q_*$ by an \textit{unknown} $\q$ and flip the sign of the term right after it. The modification of the numerical implementation is a minor one. For example, the corresponding Jacobian of this nonlinear system is replaced by
$$
D(\s,\q;\q_*) = 
\left(
\begin{array}{c|c}
I_3 & -G(\q;\q_*) \\ \hline
-H(\s,\q;\q_*) & I_3-\cos \left( \frac{h}{2} \|\s\| \right) K(\q;\q_*)
\end{array}
\right) \, .
$$
with the same matrices $G$ and $K$ as defined above, and 
\begin{eqnarray*}
H_{i,j}&=&  \frac{h}{2} \sin \left( \frac{h}{2}\|\s\| \right) \frac{s_j}{\|\s\|} \, \mbox{SLERP}_i\left(\q_*,\q,\frac{1}{2}\right)\\
& & - \frac{h}{2} \cos\left(\frac{h}{2} \|\s\|\right) \frac{s_is_j}{\|\s\|^2} - \sin \left(\frac{h}{2} \|\s\|\right) \left( \frac{\delta_{i,j}}{\|\s\|} - \frac{s_i s_j}{\|\s\|^3} \right) \, . 
\end{eqnarray*}

It might also be tempting to replace the SLERP interpolation with the simple projection, which results in the following system of nonlinear equations for both $(\s,\q)$
\begin{equation}
 \left\{
\begin{array}{l}
\s - g\left( \frac{\q_*+\q}{\|\q_*+\q\|} \right) =0 \, , \\
-\q_* + \cos\left( \frac{h}{2} \|\s\|\right) \frac{\q_*+\q}{\|\q_*+\q\|} - \sin\left(\frac{h}{2} \|\s\| \right) \frac{\s}{\|\s\|} =0  \, .
\end{array}
\right.
\label{Eqn:CNNotWork}
\end{equation}
Unfortunately, the system does not have a unique solution in $\mathbb{R}^3$ due to the projection step in both equations. The main issue is that the point $\bar{\mathbf{q}} = (\mathbf{q}_* + \mathbf{q}) / \|\mathbf{q}_* + \mathbf{q}\|$ is not unique in $\q$. To illustrate this, let us consider fixing a point $\bar{\mathbf{q}} \in \mathbb{S}^2$ and determining all $\mathbf{q} \in \mathbb{R}^3$ that will be mapped to $\bar{\mathbf{q}}$ after the averaging with $\mathbf{q}_*$ and the normalization step. This setup is depicted in Figure \ref{Fig:SphericalCN}(b). It can be observed that any points along the solid red straight line will satisfy the constraint. In particular, if we consider the two-dimensional plane (the plane containing the great circle of the points $\mathbf{q}$ and $\mathbf{q}_*$) and let $\mathbf{q}_* = (q_{*1}, q_{*2})$ and $\bar{\mathbf{q}} = (1, 0)$, it is evident that any point $(2t - q_{*1}, q_{*2})$ for all $t > 0$ will be mapped to $\bar{\mathbf{q}}$. This demonstrates the non-uniqueness of the solution.

\section{Numerical Examples}
\label{Sec:Examples}

In this section, we conduct a comparative analysis between our proposed implicit spherical integrators and several other numerical approaches. Through numerical demonstrations, we aim to demonstrate the accuracy and stability of the proposed algorithms. 

\begin{figure}[!htb]
\centering
\includegraphics[trim=0 0 10 0, clip, width=0.6\textwidth]{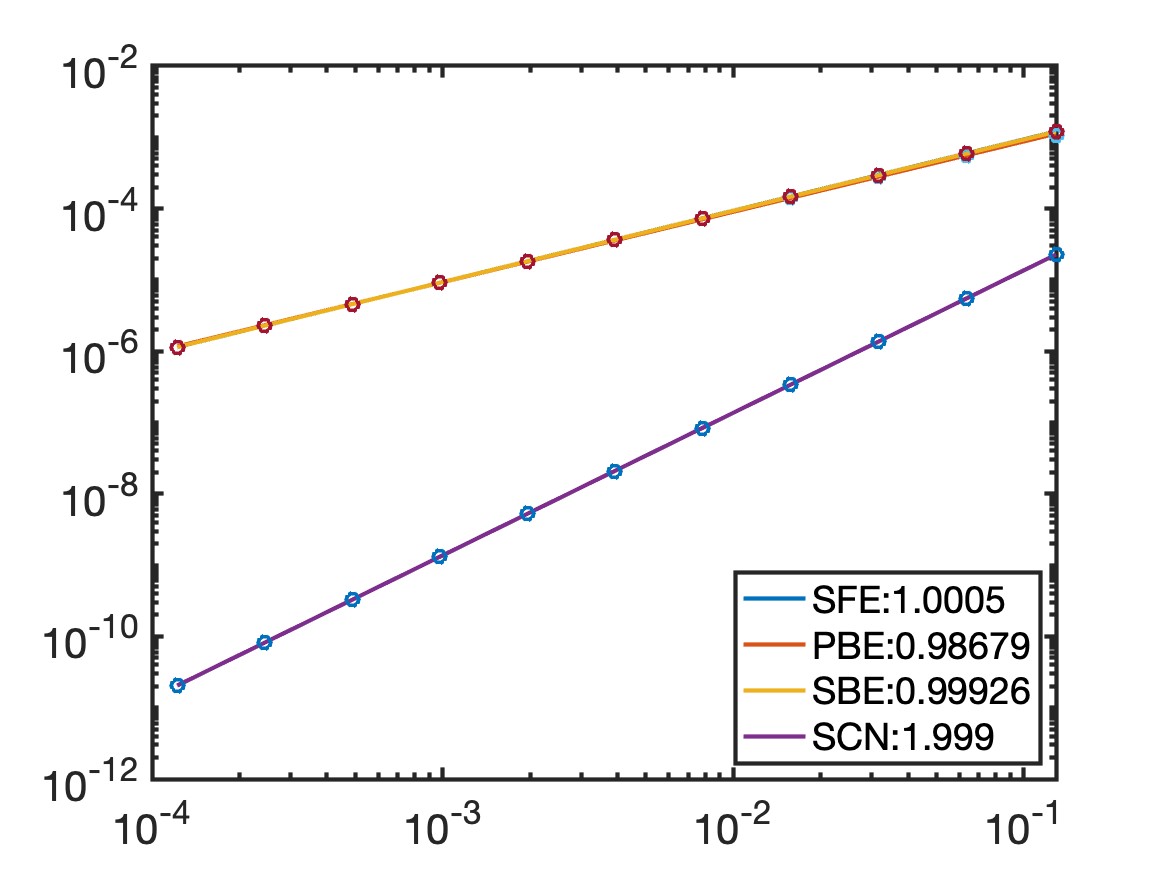}
\caption{(Section \ref{SubSec:ExConvergence}) The $E_2$ errors in the solutions obtained by Spherical-FE, Projected-BE, and Spherical-BE demonstrate first-order accuracy. In contrast, the one obtained by our proposed Spherical-CN demonstrates second-order accuracy. 
} 
\label{Fig:Convergence123}
\end{figure}

\subsection{Convergence}
\label{SubSec:ExConvergence}

This example considers the following four-point vertices flow given by
$$
f(\x)=\sum_{i=1}^4 \frac{\x_i \times \x}{2(1-\x_i\cdot \x)}
$$
for $\x_1=(1,-1,1)/\sqrt{3}$, $\x_2=(1,-1,-1)/\sqrt{3}$, $\x_3=(-2,1,0)/\sqrt{5}$, and $\x_4=(-1,-1,0)/\sqrt{2}$. The initial condition is given by $\p_0=(1,0,0)$. The exact solution is computed using the STVDRK3 scheme developed in \cite{leuchalee24} with a significantly smaller step size.

Figure \ref{Fig:Convergence123} illustrates the error in the final arrival location at $T=2$ for the solutions computed using the spherical FE method (as developed in \cite{leuchalee24}) as well as our proposed PBE, SBE, and SCN methods. We observe that all three first-order methods yield very similar solutions. The three least-squares fitted lines, with slopes close to one, almost overlap. The straight line at the bottom represents the least-squares fitted line associated with the error in the numerical solutions computed using the SCN method. We can observe that this method is second-order accurate.


\begin{figure}[!htb]
\centering
\includegraphics[width=0.99\textwidth]{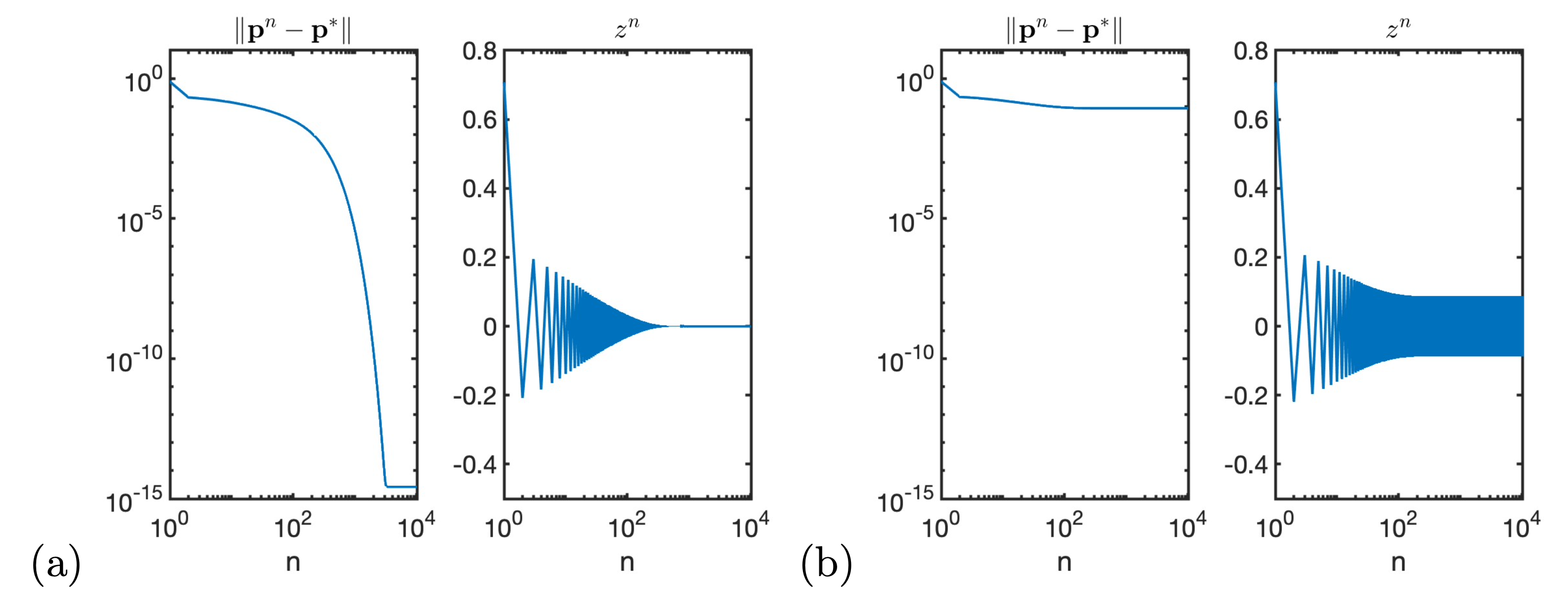}
\caption{(Section \ref{Ex:Stability} \cite{leuchalee24}) The numerical solutions were obtained using SFE with step sizes of (a) 1.99 where $z^n$ converges to 0 and (b) 2.01 where $z^n$ diverges.} 
\label{Fig:Stability1}
\end{figure}

\begin{figure}[!htb]
\centering
\includegraphics[width=0.99\textwidth]{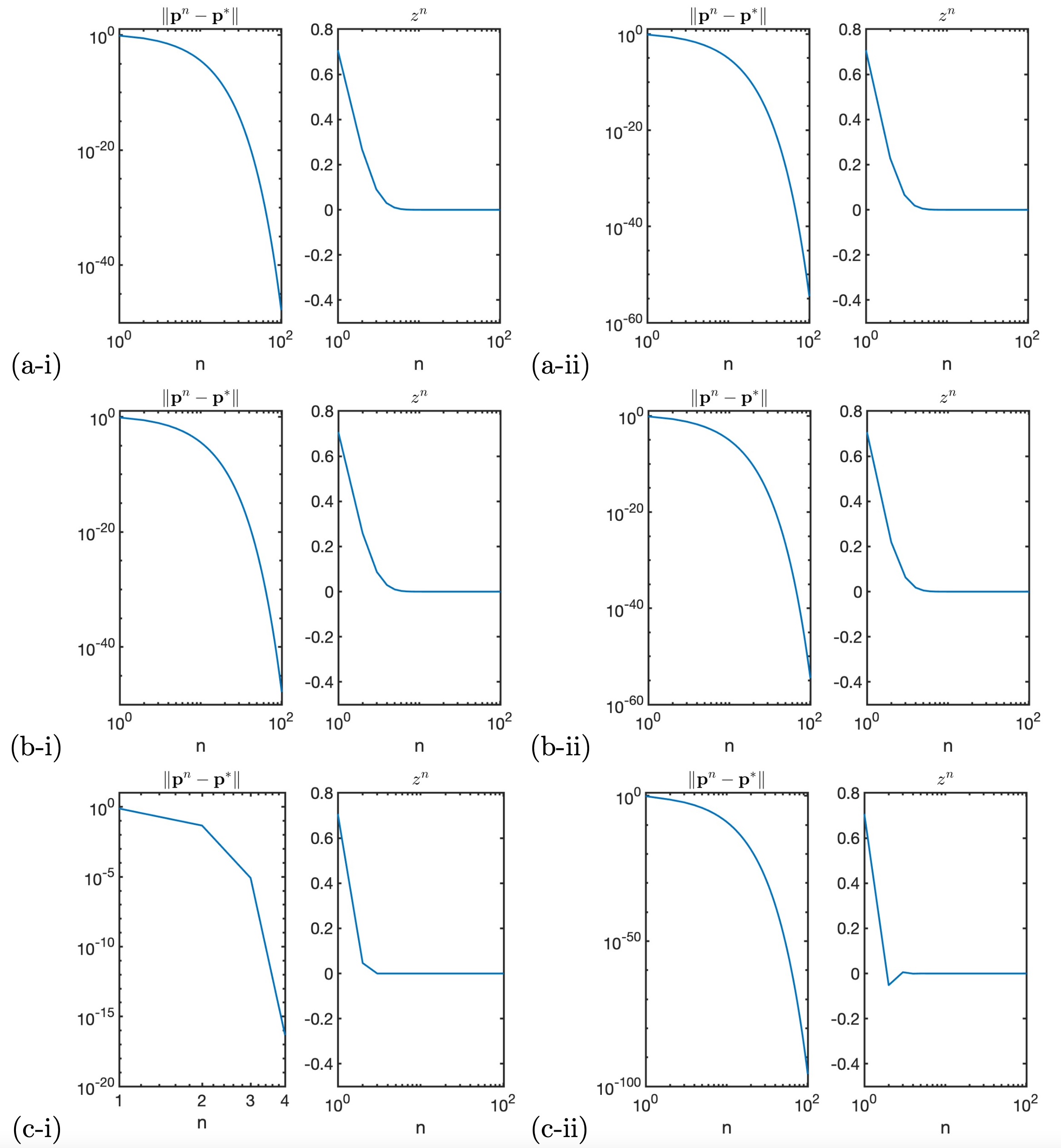}
\caption{(Section \ref{Ex:Stability}) The numerical solutions were obtained using (a) spherical backward Euler, (b) projected backward Euler, and (c) spherical Crank-Nicolson with step sizes of (i) 2.0 and (ii) 2.5. All $z^n$ converges to 0.} 
\label{Fig:Stability2}
\end{figure}

\subsection{Stability}
\label{Ex:Stability}

We demonstrate the significance of satisfying the A-stability condition in determining an appropriate time step size for a numerical method. As elaborated upon in \cite{leuchalee24}, we consider the nonlinear ODE problem $\mathbf{q}'=(I - \mathbf{q}\mathbf{q}^T ) M \mathbf{q}$ with
$$
M=\left(
\begin{array}{ccc}
\frac{1}{2} & 0 & 0 \\
0 & -\frac{1}{2} &0 \\
0 & 0 & -\frac{1}{2}
\end{array}
\right) \, .
$$
This particular model problem possesses an equilibrium point at the origin, and the matrix $M$ is characterized by three eigenvalues. Among these, 1/2 and -1/2 each have a multiplicity of 2. The eigenvector associated with the eigenvalue 1/2 aligns with the $x$-direction (i.e., $\mathbf{e}_1$), contributing to the diverging component of the solution for the linear homogeneous problem. The vectors $\mathbf{p} = \mathbf{e}_1$ and $-\mathbf{e}_1$ give rise to two equilibrium points on the unit sphere. The Jacobian matrix associated with these points on the tangent plane yields two eigenvalues of -1, indicating both equilibria are stable attractors. This highlights the importance of understanding the stability characteristics of numerical methods, especially when applied to problems with nonlinear dynamics.

Figure \ref{Fig:Stability1} presents numerical solutions obtained using the SFE method developed in \cite{leuchalee24}, with two different step sizes. One step size satisfies the A-stability condition, while the other slightly exceeds the threshold. For example, we choose a time step size $h$ less than 2 and perform simulations with $h=1.99$ and $2.01$, as shown in Figure \ref{Fig:Stability1}(a) and (b), respectively. In each figure, the left subplot shows the distance between intermediate solutions and the attractor $\mathbf{e}_1$ as a function of the iteration number, while the right subplot depicts the third component of the solution ($z_3$). The solution with $h=1.99$ demonstrates favorable convergence towards the point $\mathbf{e}_1$, while the solution with $h=2.01$ diverges. For all implicit integrators, including the SBE, PBE, and SCN methods, we examine solutions with $h=2$ and $2.5$, as illustrated in Figure \ref{Fig:Stability2}. We observe that all numerical solutions are stable for this stiff system.


\begin{figure}[!htb]
\centering
\includegraphics[width=0.99\textwidth]{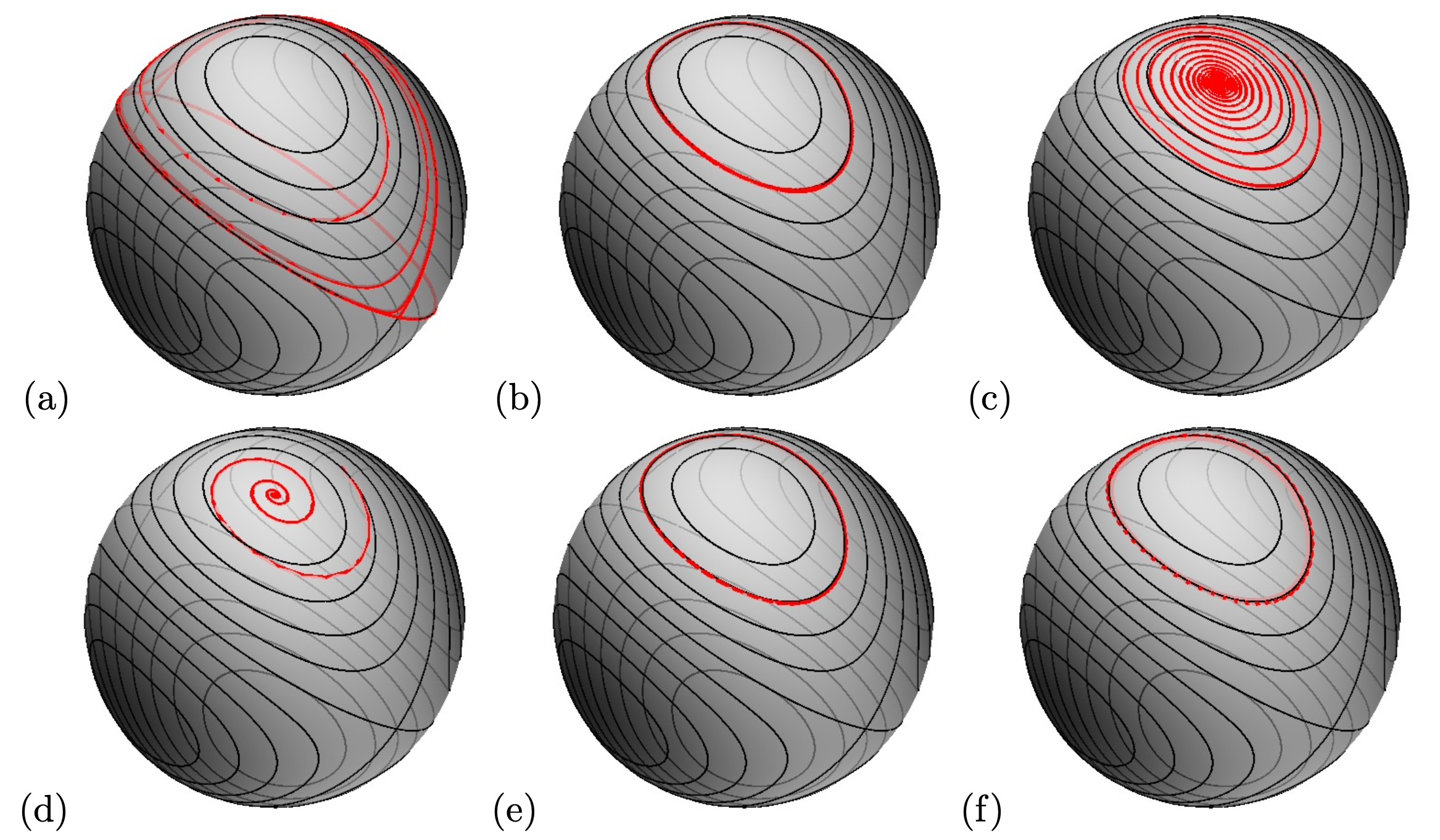}
\caption{(Section \ref{SubSec:Ham}) Numerical solution of the Hamiltonian flow obtained by two explicit schemes, (a) Spherical FE and (b) STVDRK3, with $h=0.5$ \cite{leuchalee24}. We use $h=0.5$ and solve the equations up to the final time $T=500$. 
We also show our proposed spherical backward Euler scheme using (c) $h=0.1$ and (d) $h=0.5$, and the solution using spherical Crank-Nicolson with step sizes of (e) $h=0.5$ and (f) $h=1.0$. We plot the numerical solutions in a red solid curve, while some Hamiltonian contours are in black solid lines.} 
\label{Fig:HamSolution}
\end{figure}

\begin{figure}[!htb]
\centering
\includegraphics[width=0.99\textwidth]{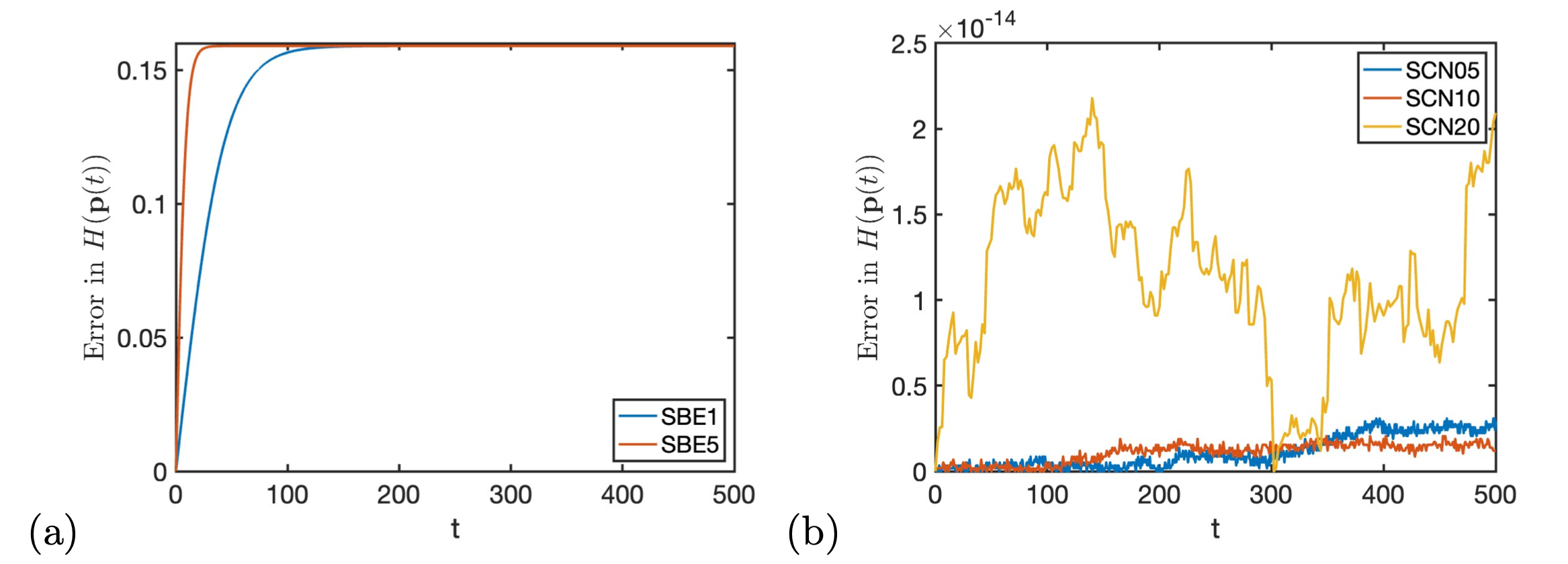}
\caption{(Section \ref{SubSec:Ham}) (a) The Hamiltonian along the solution trajectory obtained by spherical backward Euler scheme using $h=0.1$ and $h=0.5$. (c) The Hamiltonian along the solution trajectory was obtained by spherical Crack-Nicolson using $h=0.5$, $1.0$, and $h=2.0$. We solve the equations up to the final time $T=500$.} 
\label{Fig:HamSolutionHam}
\end{figure}

\begin{figure}[!htb]
\centering
\includegraphics[width=0.99\textwidth]{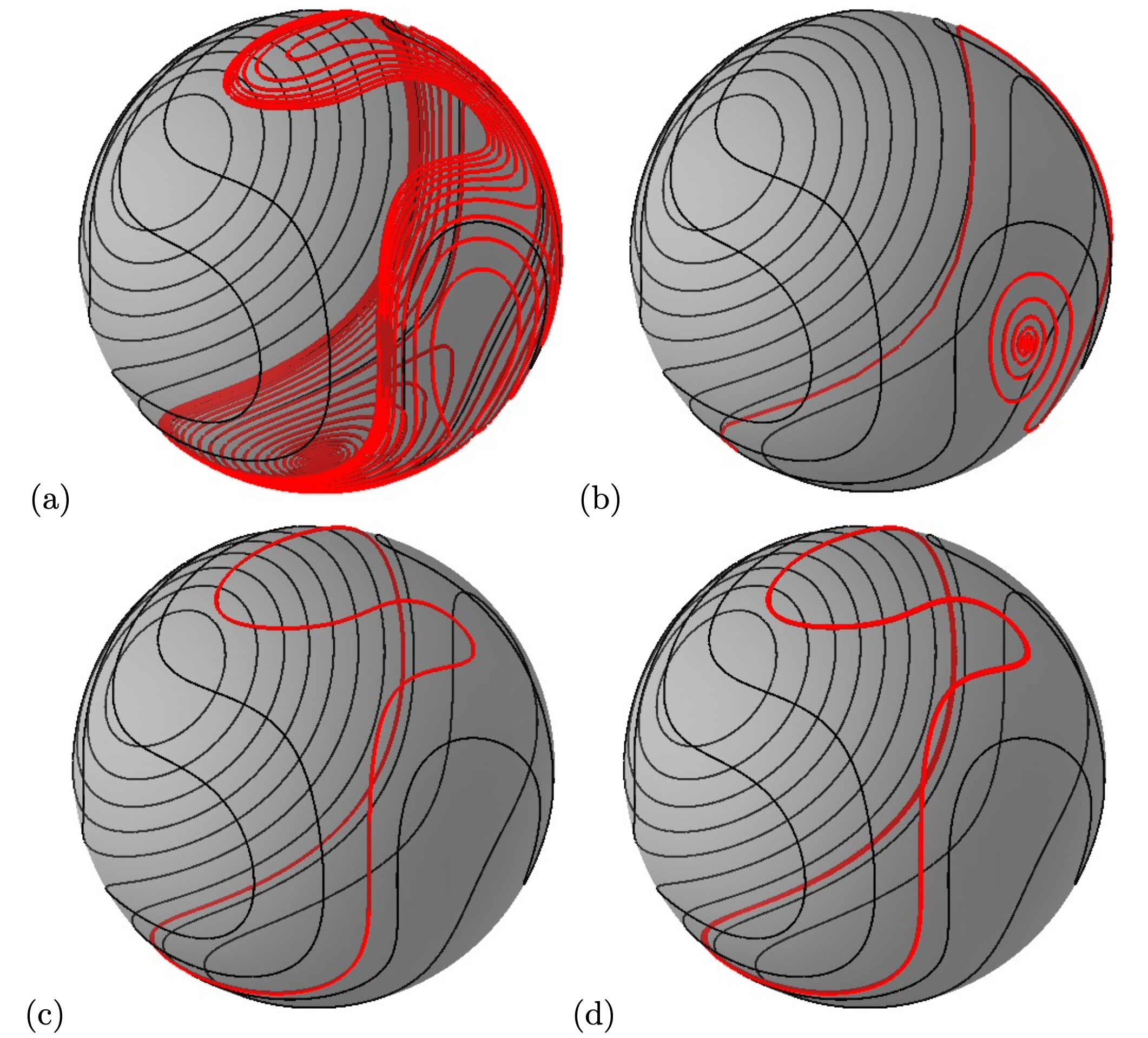}
\caption{(Section \ref{SubSec:Ham}) Numerical solution of the Hamiltonian flow obtained by our proposed spherical backward Euler scheme using (a) $h=0.1$ and (b) $h=0.5$, and the solution using spherical Crank-Nicolson with step sizes of (c) $h=0.1$ and (d) $h=0.5$. We plot the numerical solutions in a red solid curve, while some Hamiltonian contours are in black solid lines.} 
\label{Fig:HamSolution2}
\end{figure}

\begin{figure}[!htb]
\centering
\includegraphics[width=0.99\textwidth]{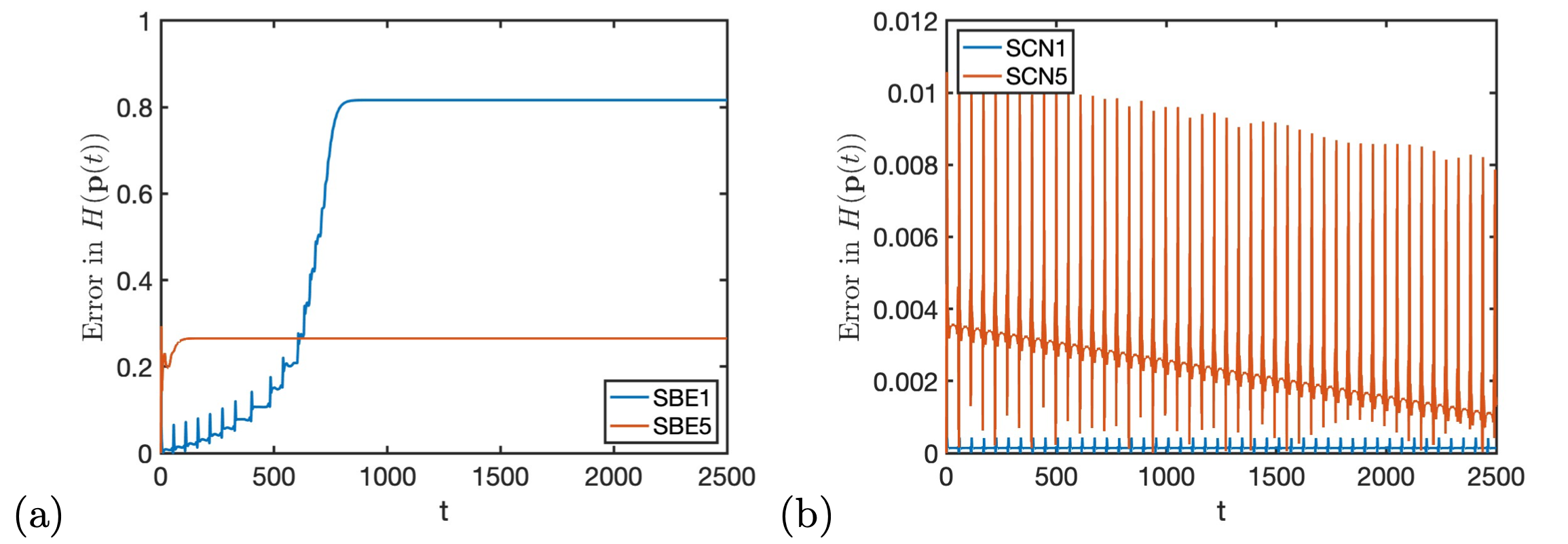}
\caption{(Section \ref{SubSec:Ham}) (a) The Hamiltonian along the solution trajectory obtained by spherical backward Euler scheme using $h=0.1$ (in blue) and $h=0.5$ (in red). (b) The Hamiltonian along the solution trajectory was obtained by spherical Crack-Nicolson using $h=0.1$ and $h=0.5$. We solve the equations up to the final time $T=2500$.} 
\label{Fig:HamSolutionHam2}
\end{figure}

\subsection{Hamiltonian Flows}
\label{SubSec:Ham}

In this example, we first consider the following Hamiltonian flow \cite{hai00}:
$$
y_1' = a_1 y_2 y_3 \, , \, y_2' = a_2 y_3 y_1 \, , \, y_3' = a_3 y_1 y_2 \, .
$$
Here, $a_1=(I_2-I_3)/(I_2 I_3)$, $a_2=(I_3-I_1)/(I_3 I_1)$, and $a_3=(I_1-I_2)/(I_1 I_2)$, where $I_1$, $I_2$, and $I_3$ are the principal moments of inertia of a rigid body. The vector $(y_1,y_2,y_3)$ physically represents the angular momentum of the rigid body frame. The above system is solved with the initial condition $(y_1(0),y_2(0),y_3(0))\in \mathbb{S}^2$. One can show that the solution to this system stays on $\mathbb{S}^2$ while the motion keeps the Hamiltonian $H(y_1,y_2,y_3)=\frac{1}{2} \left( \frac{y_1^2}{I_1} + \frac{y_2^2}{I_2} + \frac{y_3^2}{I_3} \right)$ constant. Although the exact solution conserves the Hamiltonian in the time evolution, we do not build the mechanism in the numerical scheme. We constrain only that the solution $\mathbf{p}(t)\in\mathbb{S}^2$, but not $\mathbf{p}(t)\in\mathbb{S}^2 \cap \{ H(y_1,y_2,y_3)=H(\mathbf{p}(0))\}$. This example follows \cite{hai00} and chooses $I_1=2$, $I_2=1$, and $I_3=2/3$. We also use the same initial condition $(\cos(1.1),0,\sin(1.1))$ and the timestep $h=0.5$. We solve the system up to a large final time $T=500$.

Figure \ref{Fig:HamSolution} presents the numerical solutions obtained by various methods for the Hamiltonian flow problem, including the spherical forward Euler method and the STVDRK3 approaches developed in \cite{leuchalee24}, as well as our proposed spherical backward Euler method and the spherical Crank-Nicolson scheme. Since these integrators automatically preserve the constraint, we can observe that all solutions stay on the unit sphere. However, the performance of the SBE method is unsatisfactory, as shown in Figure \ref{Fig:HamSolution}(c-d). Although these solutions have unit length, they do not preserve the Hamiltonian system well. Suppose the numerical solution preserves the system's energy, we should expect the red curve to coincide with a level contour of the Hamiltonian on the sphere, represented by a black solid line. However, we can see that both SBE solutions with different step sizes converge to a single point on the sphere. On the other hand, when using the same step size (i.e., $h=0.5$), the SCN solution exhibits behavior similar to the STVDRK3 (as shown in Figure \ref{Fig:HamSolution}(b)), which is a third-order scheme. Additionally, the solution with an even larger step size still appears stable, as shown in Figure \ref{Fig:HamSolution}(f).

To further investigate energy conservation, we plot the Hamiltonian as a function of time in Figure \ref{Fig:HamSolutionHam}. In Figure \ref{Fig:HamSolutionHam}(a), we show the relative error in the Hamiltonian along the trajectories computed by the SBE method as a function of time. These solutions exhibit more than a 15\% error as they converge to a single point on the unit sphere. In contrast, the symmetric integrator SCN produces Hamiltonian-preserving numerical solutions for various step sizes. In Figure \ref{Fig:HamSolutionHam}(b), we plot the relative error in the Hamiltonian for step sizes of $h=0.5$, $1.0$, and $2.0$. All errors are on the order of machine epsilon, which indicates that the energy of the system is well-preserved. This suggests that developing a symplectic integrator for general Hamiltonian systems may be possible based on this symmetric integrator.

We have also considered the flow induced by $H(\mathbf{p}) = \frac{1}{2} \sum_{j=1}^3 \frac{1}{I_j} \left( p_j^2 + \frac{2}{3} p_j^3 \right)$, where $I_1=1$, $I_2=2$, and $I_3=4$ \cite{mclmodver14,ceor20}, representing a nonlinear perturbation of a spinning top. In Figure \ref{Fig:HamSolution2}, we present the solutions computed using our proposed SBE and SCN methods with different step sizes up to a final time $T=2500$. Similar to the previous case, we observe that the solution obtained by the first-order scheme does not form a closed trajectory on the sphere. These solutions appear to be attracted to one of the equilibria of the system. In contrast, the solutions obtained using the SCN method are periodic, as shown in Figure \ref{Fig:HamSolution2} (c-d). The solutions obtained using the symmetric integrator exhibit closed trajectories on the unit sphere, showcasing the effectiveness of the method in preserving the geometric properties of the Hamiltonian system. It is worth noting that these geometric properties are not explicitly imposed in the design of the algorithm, further highlighting the robustness and intrinsic preservation capabilities of the symmetric integrator. Figure \ref{Fig:HamSolutionHam2} presents the relative error in the Hamiltonian of the computed solutions as a function of time. We observe that the integrator preserves the energy of the system much better compared to the first-order integrator. This result is consistent with the example of the rigid body rotation system, indicating the superior energy preservation properties of the seemingly symplectic integrator.

\section*{Acknowledgment}
The work was supported by the Hong Kong RGC under grants 16302223 and 16300524.

\section*{Declarations}

\paragraph{Availability of data and materials} Not applicable.

\paragraph{Conflict of interest} The authors declare that they have no conflict of interest.

\bibliographystyle{plain}
\bibliography{syleung}

\end{document}